\documentclass[12pt]{amsart}
\usepackage{amsfonts,amssymb,amsmath,mathrsfs}

 2   
\font\Bbb=msbm10 scaled 1200

\RequirePackage{ifpdf}
\RequirePackage{xr-hyper}
\RequirePackage{hyperref}

\newtheorem{thm}{Theorem}[section]

\newtheorem{prop}[thm]{Proposition}
\newtheorem{cor}[thm]{Corollary}

\newtheorem{rmk}[thm]{Remark}

\newcommand{\thmref}[1]{Theorem~\ref{#1}}

\newcommand{\rmkref}[1]{Remark~\ref{#1}}
\newcommand{\propref}[1]{Proposition~\ref{#1}}

\begin{document}

\title
[L-functions of half-integral weight]
{Non-vanishing of $L$-functions associated to cusp forms of 
half-integral weight} 
\author[B. Ramakrishnan and K. D. Shankhadhar]{B. Ramakrishnan and Karam Deo Shankhadhar}
\address[B. Ramakrishnan]{Harish-Chandra Research Institute, 
Chhatnag Road, Jhunsi, Allahabad 211 019 (India)}
\email{ramki@hri.res.in}
\address[Karam Deo Shankhadhar] {Harish-Chandra Research Institute, 
Chhatnag Road, Jhunsi, Allahabad 211 019 (India)\\
{\sl Current Address:} The Institute of Mathematical Sciences, IV cross road, CIT Campus, Taramani, Chennai 600 113 (India)}
\email{karam@imsc.res.in}


\maketitle

\begin{abstract}
In this article, we prove non-vanishing results for $L$-functions associated to holomorphic cusp forms of 
half-integral weight on average (over an orthogonal basis of Hecke eigenforms). This extends a result of 
W. Kohnen \cite{kohnen} to forms of half-integral weight. 
\end{abstract}

\section{Introduction} 

In \cite{kohnen}, W. Kohnen proved that, given a real number $t_0$ and 
a positive real number $\epsilon$, for all $k$ large enough the sum of the 
functions $L^*(f,s)$ 
with $f$ running over a basis of (properly normalized) 
Hecke eigenforms of weight $k$ does not vanish
on the line segment Im $s=t_0,(k-1)/2<{\rm{Re}}(s)<k/2-\epsilon,k/2+\epsilon
<{\rm{Re}}(s)<(k+1)/2$. As a consequence, he proved that 
for any such point $s$, for $k$ large enough
there  exists a Hecke eigenform of weight $k$ on 
$SL_2({\mathbb Z})$ such that 
the corresponding $L$-function value at $s$ is non-zero.
Using similar methods, in \cite{raghuram}, A. Raghuram generalised Kohnen's 
result for the average of $L$-functions over a basis of newforms 
(of integral weight) of level $N$ with primitive character modulo $N$. 
In this article, we extend Kohnen's method to forms of 
half-integral weight. As a consequence, we show that for any given point $s$ inside the critical strip ~$k/2-1/4<{\rm Re}(s)< k/2+3/4$,  
there  exists a Hecke eigen cusp form $f$ of half-integral weight $k+1/2$ on $\Gamma_0(4N)$ with 
character $\psi$  such that the corresponding $L$-function value at $s$ is non-zero, and the first Fourier 
coefficient of $f$ is non-zero. It should be noted that the normalisation of Fourier coefficients 
of forms of half-integral weight is still an open question. Also, contrary to the result of Kohnen in the case of integral weight 
modular forms on $SL_2({\mathbb Z})$, we get the non-vanishing result inside the critical strip including the central line 
(see \rmkref{rmk:00}). Our results are obtained for $N$  sufficiently large if $k$ is fixed and vice versa. 
In particular when $N=1$, for sufficiently large $k$, either $f$ is a newform in the full space 
or $f$ is a Hecke eigenform in the Kohnen plus space.
 
\section{Notations and Main Theorems}

Let $N\ge 1$, $k\ge 3$ be  integers and $\psi$ be an even Dirichlet 
character modulo $4N$.  
Let $S_{k+1/2}(4N,\psi)$ be the space of cusp forms of weight 
$k+1/2$, level $4N$ with character $\psi$ (\cite{koblitz}, \cite{shimura}). 
Let $L(f,s)$ be the $L$-function 
associated to the cusp form $f \in S_{k+1/2}(4N,\psi)$ defined by 
$L(f,s) = \sum_{n\ge 1} a_f(n) n^{-s}$,
where $a_f(n)$ denotes the $n$-th Fourier coefficient of $f$. Then by 
\cite[Proposition 1]{mmr}, the completed $L$-function defined by 
$L^*(f,s) := (2\pi)^{-s}(\sqrt{4N})^s\Gamma(s) L(f,s)$ 
has the following functional 
equation
\begin{equation}\label{eq:1}
L^*(f\vert H_{4N}, k+1/2 -s) = L^*(f,s),
\end{equation}
where $H_{4N}$ is the Fricke involution on $S_{k+1/2}(4N,\psi)$ defined by
$$
f|H_{4N}(z) = i^{k+1/2}(4N)^{-k/2-1/4}z^{-k-1/2}f(-1/4Nz).
$$
For $f, g \in S_{k+1/2}(4N,\psi)$, let $\langle f,g\rangle$ denote the 
Petersson scalar product of $f$ and $g$. It is known 
that the  space $S_{k+1/2}(4N,\psi)$ has an orthogonal basis of 
Hecke eigenforms with respect to all Hecke operators 
$T(p^2)$, $p\not\vert 2N$. Let $\{f_1, f_2, \ldots, f_d\}$ be such an 
orthogonal basis of Hecke eigenforms, where $d$ is the dimension of the space $S_{k+1/2}(4N,\psi)$ (see for 
example \cite{shimura}).  Let $K$ be the operator defined  by  $f\vert K(z) = \overline{f(-\overline{z})}$.   
Since $K H_{4N} = H_{4N} K$ on 
$S_{k+1/2}(4N,\psi)$, we have $f\vert (K H_{4N})^2 = f$. Also, the operators $K$ and 
$H_{4N}$ commute with the Hecke operators $T(p^2)$, $p\not\vert 2N$. Therefore, 
for the basis elements $f_j, ~1\le j\le d$, we may assume that $f_j\vert K H_{4N} 
= \lambda_{f_j} f_j$, where $\lambda_{f_j} = \pm 1$. 

\smallskip

We now state the main results of this article.

\begin{thm}\label{thm1}
Let $N\ge 1$ be a fixed integer. 
Let $\{f_1,f_2,\ldots,f_d\}$ be an orthogonal basis 
as above. Let $r_0\in {\mathbb R}$.
Then there exists a constant 
$C = C(r_0)$ depending only on $r_0$  such that for $k> C$, 
the function   
\begin{equation*}
\sum_{j=1}^{d}\frac{L^*(f_j, s)}{\langle f_j,f_j\rangle} \lambda_{f_j} a_{f_j}(1)
\end{equation*}
doesn't vanish for any point $s=\sigma+i r_0$ with $k/2-1/4<\sigma<k/2+3/4$.
\end{thm}

\begin{thm}\label{thm2}
Let $k\ge 3$ be a fixed integer. 
Let $\{f_1,f_2,\ldots,f_d\}$ be an orthogonal basis 
as above. Let $r_0\in {\mathbb R}$. Then there exists a constant 
$C' = C'(r_0)$ depending only on $r_0$ such that for $N> C'$, 
the function 
\begin{equation*}
\sum_{j=1}^{d}\frac{L^*(f_j, s)}{\langle f_j,f_j\rangle} \lambda_{f_j} a_{f_j}(1)
\end{equation*}
doesn't vanish for any point $s=\sigma+i r_0$ 
with $k/2-1/4<\sigma<k/2+3/4$. 
\end{thm}

The following corollary is an easy consequence of the above 
two theorems.

\begin{cor}\label{cor:1}
Let $s_0$ be a point inside the critical strip $k/2-1/4 <{\rm Re}(s_0)<k/2+3/4$. 
If either $k$ or $N$ is suitably large then there exists a Hecke eigenform $f$ belonging 
to $S_{k+1/2}(4N,\psi)$ such that $L(f,s_0) \not =0$ and $a_f(1)\not=0$.  
\end{cor}

\begin{rmk}\label{rmk:0}
{\rm Though we consider Hecke eigenforms of half-integral weight, the $L$-function corresponding to such a Hecke eigenform does not have 
an Euler product. } 
\end{rmk}

\section{Proof}
The proof is on the same lines as that of Kohnen \cite{kohnen} 
and so we give only a sketch. 
First, let us recall the Poincar\'e series in 
$S_{k+1/2}(4N,\psi)$. We define the $n$-th Poincar\'e series in $S_{k+1/2}(4N,\psi)$ 
as follows.
\begin{equation}\label{eq:2}
P_{n,k+1/2,4N,\psi}(z)=\frac{1}{2} \sum_{(c,d)\in \mathbb{Z}^2\atop{(c,d)=1
,4N|c}}\overline{\psi}(d)\left(\frac{c}{d}\right)\left(\frac{-4}{d}\right)^
{k+1/2}(cz+d)^{-(k+1/2)}e\left(n\frac{a_0z+b_0}{cz+d}\right),
\end{equation} 
where in the summation above, for each coprime pair $(c,d)$ and $4N|c$, 
we make a fixed choice of $(a_0,b_0) \in \mathbb{Z}^2$ with $a_0d-b_0c =1$ and
$e(x)$ stands for $e^{2\pi i x}$.
We have the following characterization of the Poincar\'e series.
\begin{equation}\label{eq:3}
\langle f,P_{n,k+1/2,4N,\psi} \rangle = \frac{\Gamma(k-1/2)}{i_{4N} (4\pi n)
^{k-1/2}} ~a_f(n), \qquad f \in S_{k+1/2}(4N,\psi),
\end{equation}
where $i_{4N}$ is the index of $\Gamma_0(4N)$ in $SL_2({\mathbb Z})$. 

\noindent Next, let us define the kernel function for the special values of the $L$-function 
associated to a cusp form of half-integral weight. A similar function for forms of integral 
weight was considered by Kohnen \cite{kohnen}. 
Let $z\in {\mathcal H}$ and $s\in{\mathbb C}$ with $1< \sigma <
k -1/2$, $\sigma = {\rm Re}(s)$. Define
\begin{equation}\label{eq:4}
R_{s;k,N,\psi}(z) = \gamma_k(s) {\sum}'\overline{\psi}(d)
\left(\frac{c}{d}\right) \left(\frac{-4}{d}\right)^{k+1/2}
(cz+d)^{-(k+1/2)} \left(\frac{az+b}{cz+d}\right)^{-s},
\end{equation}
where
\begin{equation}\label{eq:5}
\gamma_{k}(s) = \frac{1}{2} ~e^{\pi is/2}~\Gamma(s)\Gamma(k+1/2-s),
\end{equation}
and the sum $\sum'$ varies over all matrices $\begin{pmatrix} a&b\\
c&d \\\end{pmatrix}$ $\in \Gamma_0(4N)$.
The condition $1< \sigma < k-1/2$
ensures that the above series converges
absolutely and uniformly on compact subsets of ${\mathcal H}$ and hence it
represents an analytic function on ${\mathcal H}$. The function
$R_{s,\psi}(z) :=  R_{s; k,N,\psi}(z)$ is a cusp form in
$S_{k+1/2}(4N,\psi)$. 

\noindent For a given  $c,d \in {\mathbb Z}$ 
with $\gcd(c,d)=1$ and $4N \vert c$,
we choose $a_0, b_0$ such that $a_0d -b_0c =1$. Then any other solution
$a,b$ of $ad-bc =1$ is given by
$a = a_0+nc$ and $b = b_0 +nd$ for some $n\in {\mathbb Z}$.
Hence,
\begin{equation}
R_{s, \psi}(z) = \gamma_k(s)
\sum_{(c,d)=1\atop{4N\mid c}} \sum_{n\in {\mathbb Z}}\overline{\psi}(d)
\left(\frac{c}{d}\right) \left(\frac{-4}{d}\right)^{k+1/2}
(cz+d)^{-k-1/2} \left(\frac{a_0z+b_0}{cz+d}+n\right)^{-s}.
\end{equation}
Using Lipschitz's formula
\begin{equation}
\sum_{n=-\infty}^{\infty} (z + n)^{-s} = \frac{e^{-\pi i s/2}
(2\pi)^s}{\Gamma(s)} \sum_{n\ge 1} n^{s-1} e(nz) \qquad
(z\in {\mathcal H}, \sigma>1),
\end{equation}
we get
\begin{equation}
\begin{split}
R_{s,\psi}(z) & = \gamma_k(s)
\sum_{(c,d)=1\atop{4N\mid c}} \overline{\psi}(d)
\left(\frac{c}{d}\right) \left(\frac{-4}{d}\right)^{k+1/2}
(cz+d)^{-k-1/2} \frac{e^{-\pi is/2} (2\pi)^s}{\Gamma(s)}  \\
&  \hskip 7cm \times\sum_{n\ge 1} n^{s-1}
e\left(n\frac{a_0z+b_0}{cz+d}\right) \\
\end{split}
\end{equation}
\begin{equation*}
\begin{split}
&= (2\pi)^s \Gamma(k+1/2-s)\sum_{n\ge 1} n^{s-1} \\
&\hskip 4cm \times \frac{1}{2}
\sum_{(c,d)=1\atop{4N\mid c}} \overline{\psi}(d)
\left(\frac{c}{d}\right) \left(\frac{-4}{d}\right)^{k+1/2}~
(cz+d)^{-k-1/2} e\left(n\frac{a_0z+b_0}{cz+d}\right).
\end{split}
\end{equation*}
Here, we have used the absolute convergence of the above sum in the region $1<\sigma<k-\beta-1/2$ 
and so the interchange of summations is allowed, where $\beta = k/2-1/28$ is the exponent of the estimate for the 
Fourier coefficients of a cusp form of weight $k+1/2$ on $\Gamma_0(4N)$ obtained by H. Iwaniec \cite{iwaniec}.
Thus, for $1 < \sigma < k-\beta-1/2$,
\begin{equation}\label{char-r-p}
R_{s,\psi}(z) = (2\pi)^s \Gamma(k+1/2 - s)\sum_{n\ge 1} n^{s-1}
P_{n,k+1/2,4N,\psi}(z).
\end{equation}
Using equation \eqref{eq:3} with the last equation, we get 
for $1 < \sigma < k-\beta-1/2$,
\begin{equation*}
\langle f, R_{\overline{s},\psi}\rangle = \frac{\pi \Gamma(k-1/2)}{i_{4N} 2^{k-3/2} (4N)^{k/2+1/4-s/2}} L^*(f, k+1/2-s), 
\end{equation*}
for all $f\in S_{k+1/2}(4N,\psi)$. Using this, we have 
\begin{equation}\label{petersson}
R_{s,\psi} =  \frac{2^{-2k+1+s} \pi \Gamma(k-1/2)}{i_{4N} N^{k/2+1/4-s/2}} \sum_{j=1}^{d} 
\frac{L^*(f_j\vert K, k+1/2-s)}{\langle f_j,f_j \rangle} f_j,
\end{equation}
where the sum varies over the orthogonal basis $\{f_j\}$. 
Using $f_j\vert K H_{4N} = 
\lambda_{f_j} f_j$ together with the functional equation 
\eqref{eq:1}, 
we get for $1<\sigma<k-\beta-1/2$,
\begin{equation}\label{identity}
R_{s,\psi} =  \frac{2^{-2k+1+s} \pi \Gamma(k-1/2)}{i_{4N} N^{k/2+1/4-s/2}} \sum_{j=1}^{d} 
\frac{L^*(f_j, s)}{\langle f_j,f_j \rangle} \lambda_{f_j} f_j.
\end{equation}
This equality has been established for $1<\sigma<k-\beta-1/2$. Since the right-hand side is an entire function 
of $s$, this gives an analytic continuation of the kernel function $R_{s,\psi}$ for all $s\in {\mathbb C}$. 

\noindent
Next, we need the Fourier expansion of the function $R_{s,\psi}$. In an earlier version of \cite{mmr}, the Fourier expansion of the function 
$R_{s,\psi}(z)$ was derived, which we present here. 
Let 
$$
R_{s,\psi}(z) = \sum_{n\ge 1} a_{s,\psi}(n) e^{2\pi inz}
$$
be the Fourier expansion of $R_{s,\psi}$, where the Fourier coefficients $a_{s,\psi}(n)$ are 
given by 
\begin{equation}\label{rschi-fourier}
\begin{split}
a_{s,\psi}(n) & = (2\pi)^s \Gamma(k+1/2 -s) n^{s-1} +  e^{\pi is/2} 
(-2\pi i)^{k+1/2}  n^{k-1/2} \qquad \qquad \\
&\qquad \times \!\!\!\!\!\sum_{(a,c)\in {\mathbb Z}^2, ac\not= 0 \atop{\gcd(a,c)=1, 
4N\vert c}}
\!\!\!\!\!\!\!\!\!
{\psi}(a)\left(\frac{c}{a}\right)\left(\frac{-4}{a}\right)^{k+1/2} c^{s-k-1/2} 
a^{-s} e^{2\pi ina'/c} {}_1f_1(s,k+1/2;-2\pi in/ac),
\end{split}
\end{equation}
where $a'$ is an integer which is the inverse of $a$ modulo $c$ and 
\begin{equation}\label{1f1}
{}_1f_1(\alpha,\beta;z) = \frac{\Gamma(\alpha)\Gamma(\beta-\alpha)}{\Gamma(\beta)} {}_1F_1(\alpha,\beta;z).
\end{equation}
Here ${}_1F_1(\alpha,\beta;z)$ is the Kummer's degenerate hypergeometric function.

\noindent 
We now give the proof of our theorems.  
Assume that $\sum_{j} \frac{L^*(f_j,s)}{\langle f_j, f_j\rangle} \lambda_{f_j} 
a_{f_j}(1) = 0$, for $s$ as in the theorem. This implies that (for the 
values of $s$)  
the first Fourier coefficient of $R_{s,\psi}$ is zero. Dividing by 
$(2\pi)^s \Gamma(k+1/2 -s)$, we obtain, 
\begin{equation}\label{eq:10}
\begin{split}
& 1  +  e^{\pi is/2} \frac{(2\pi)^{k+1/2 -s}}{i^{k+1/2}\Gamma(k+1/2 -s)}
\sum_{a,c\in {\Bbb Z}, ac\not =0 \atop{(a,c)=1,4N\vert c}}\psi(a) 
\left(\frac{c}{a}\right) 
\left(\frac{-4}{a}\right)^{k+1/2} \frac{c^{s-k-1/2}}{a^{s}}\\
&\hskip 7cm \times e^{2\pi ia'/c} {}_1f_1(s,k+1/2;-2\pi i/ac) = 0.\\
\end{split}
\end{equation}
In particular, let $s= k/2+1/4-\delta+ir_0$, where $0 \leq\delta<1/2$. 
Now, one has 
\begin{equation*}
\vert {}_1f_1(s,k+1/2;-2\pi i/ac)\vert \le 1
\end{equation*}
(see \cite{kohnen}). Taking the absolute value in  \eqref{eq:10} and using the above estimate,  we get 
\begin{equation*}\label{eq:11}
\begin{split}
1 &\le A(r_0)  ~\frac{\pi^{k/2+1/4+\delta}}{\vert \Gamma(k/2+1/4+\delta- ir_0)\vert}\frac{1}
{(2N)^{k/2+1/4+\delta}} \left(\sum_{a,c\in {\mathbb Z}, ac\not=0,\atop{(a,2Nc)=1}}
\frac{1}{a^{k/2+1/4-\delta} \cdot c^{k/2+1/4+\delta}}\right)\\
&\le  A(r_0)~B  ~\frac{\pi^{k/2+1/4+\delta}}{\vert \Gamma(k/2+1/4+\delta- ir_0)\vert}\frac{1}
{(2N)^{k/2+1/4+\delta}}, 
\end{split} 
\end{equation*}
where $A(r_0)$ is a constant depending only on $r_0$ and $B>0$ is an 
absolute constant.  To prove \thmref{thm1}, we fix $N$ and allow $k$ tend to infinity and for 
the proof of \thmref{thm2}, we fix $k$ and allow $N$ tend to infinity. In either case, the right-hand 
side goes to zero (for fixed $N$ one should use the Stirling's approximation),  a contradiction. 
This completes the proof. 


\section{Remarks}\label{remark}

\begin{rmk}\label{rmk:00}
{\rm In \cite{kohnen}, the non-vanishing result was obtained for $s$ inside the critical strip with the condition that ${\rm Re}(s) \not= k/2$, the center 
of the critical strip. However, since the level of the modular forms considered in this paper is greater than $1$, we need not assume this condition and our 
results are valid for all $s$ inside the critical strip ~$k/2-1/4<{\rm Re}(s)<k/2+3/4$. 
We also remark that the same is true in \cite{raghuram}, since the level $M$ is greater 
than $1$. The reason is as follows. When $M=1$ (i.e., when one considers the case of forms of integral weight on $SL_2({\mathbb Z})$), 
while deriving the Fourier expansion of the function $R_{s}$, the term corresponding to $ac=0$ has two contributions ($c=0$ and $a=0$). 
Therefore, in the estimation of the first Fourier coefficient, there is an extra term on the right-hand side (see \cite[p. 189, Eq.(10)]{kohnen}).  
Due to the appearane of this extra term, in order to get a contradiction, the central values have to be omitted. 
Since the case $a=0$ doesn't arise for the levels $M>1$, we do not get the extra term on the right-hand side in the estimation. Therefore, this gives the 
advantage of considering all the values of $s$ inside the critical strip. In particular, one obtains non-vanishing results for forms at the 
center of the critical strip when the level is greater than $1$.
}
\end{rmk}

\begin{rmk}\label{rmk:epsilon}
{\rm 
The average sum in \thmref{thm1} (and also in \thmref{thm2}) contains an extra factor $\lambda_{f}$ (and of course the first Fourier coefficient), 
which does not appear in Kohnen's result. In the case of level $1$, we have the functional equation $L^*(f,k-s) = (-1)^{k/2}L^*(f, s)$, whereas 
when the level $M$ is greater than $1$, we have a different functional equation in the sense that on the one side we have $L^*(f,s)$ and on the 
other side we have $L^*(f\vert H_M, k-s)$ and therefore, in the final form of the functional equation, the root number will depend on the function, 
especially the eigenvalue of $f$ under the Fricke involution $H_M$. (In the case of half-integral weight, the Fricke involution 
is $H_M=H_{4N}$ in our notation.) Therefore, we will have an extra factor, which we call $\lambda_{f}$. 
Note that the extra factor corresponding to the eigenvalues under $H_M$ also appears in Raghuram's 
results (see \cite{raghuram}).  Since normalization of Fourier coefficients is not known in the case of half-integral weight, we also have the
first Fourier coefficients appearing in the average sum.
}
\end{rmk}

\begin{rmk}\label{rmk:1}
{\rm 
Let us consider the space $S_{k+1/2}(4N,\psi)$, where $\psi$ is an even primitive Dirichlet 
character modulo $4N$.  Then it is known from the work of Serre 
and Stark \cite{se-st} that the space $S_{k+1/2}(4N,\psi)$ is the space of newforms. 
Hence, the orthogonal basis consists of newforms. 
In this case, the Hecke eigenform $f$ in Corollary \ref{cor:1} will be a newform of level $4N$. 
}
\end{rmk}

\begin{rmk}\label{rmk:2}
{\rm 
Let us consider the case $N=1$. That is, we consider the space $S_{k+1/2}(4)$. 
Let $\{f_1,f_2,\ldots, f_{d_1}\}$ be an orthogonal  basis of $S_{k+1/2}^{new}(4)$ which are newforms 
(see \cite{mrv}).    Let $\{g_1, g_2, \ldots, g_{d_2}\}$ be an orthogonal basis of $S_{k+1/2}^+(4)$ (see \cite{kohnen-plus}), which are 
Hecke eigenforms such that the set $\left\{g_1\pm g_1\vert W(4),  \ldots,  g_{d_2}\pm g_{d_2}\vert W(4)\right\}$ forms an orthogonal 
basis of $S_{k+1/2}^{old}(4)$. Here,  $d_1+2d_2 =d$ is the dimension of the space $S_{k+1/2}(4)$ and 
$W(4)$ is the Atkin-Lehner  $W$-operator for the prime $p=2$ on $S_{k+1/2}(4)$. 
Thus, an orthogonal basis of Hecke eigenforms for the space $S_{k+1/2}(4)$ is given as follows:
$$
\left\{f_1, f_2, \ldots, f_{d_1}, g_1\pm g_1\vert W(4), g_2\pm g_2\vert W(4), \ldots, g_{d_2}\pm g_{d_2}\vert W(4)\right\}. 
$$
In this case, we get the following result as a consequence 
of \thmref{thm1}.  
Let  $s_0$ be a point inside the critical strip $k/2-1/4 <{\rm Re}(s_0)<k/2+3/4$. 
If $k$ is suitably large, then there exists a $j$, with $1\le j\le d_1$ or $1\le j\le d_2$ such that 
\begin{equation}\label{main}
L(f_j,s_0)\not=0, a_{f_j}(1)\not=0 \quad {\rm or}\quad L(g_j\pm g_j\vert W(4),s_0)\not=0, a_{g_j}(1)\pm 2^{-k} a_{g_j}(4)\not=0.
\end{equation}
For the last assertion in the above equation \eqref{main}, we use the fact  that $W(4) = 2^{-k} U(4)$ on $S_{k+1/2}^+(4)$, where 
$U(4)$ is the Hecke operator for $p=2$ on $S_{k+1/2}(4)$. Note that $W(4) = H_4$ on $S_{k+1/2}(4)$ and therefore, in the second 
case of \eqref{main}, for a $j$ with $1\le j\le d_2$, it follows from the functional equation that either $L(g_j,s_0)\not=0$ or 
$L(g_j, k+1/2-s_0)\not=0$. Hence, for any given point $s$ inside the critical strip, our theorem gives (for sufficiently large $k$) 
the existence of a newform $f$ in  $S_{k+1/2}^{new}(4)$ such that $L(f,s)\not=0$ 
or a Hecke eigenform $g$ in the plus space $S_{k+1/2}^+(4)$ such that $L(g,s)\not=0$ or $L(g,k+1/2-s)\not=0$.  
Correspondingly, we also get the non-vanishing of the first Fourier coefficient (if it is a newform in $S_{k+1/2}^{new}(4)$) or the 
first or the $4$-th Fourier coefficient (if it is a Hecke eigenform in $S_{k+1/2}^+(4)$). 
}
\end{rmk}

\smallskip

\subsection{Additional remarks:} 

In view of a recent result of N. Kumar \cite{NK}, we make the following additional remarks. 

Let $N$ be any positive integer and let ~$E(k+1/2,4N)$~ denote the set of all
Hecke eigenforms ~$h\in S_{k+1/2}(4N)$~ 
(the vector space of cusp forms of weight $k+1/2$ on $\Gamma_0(4N)$ with trivial character) such that the ~$L$-value ~$L(h, k/2+1/4)\neq 0$.
In \cite[\S 4]{NK}, the following results are obtained. 
\begin{thm}\label{thm3}
For ~$h\in E(k+1/2,4N)$~ one has
\begin{equation*}
\frac{L'(h,k/2+1/4)}{L(h,{k}/{2}+{1}/{4})}=-\Psi(k/2+1/4)+\log(\pi),
\end{equation*}
where ~$\Psi$~ is the logarithmic derivative of the gamma function ~$\Gamma$. Further, for such an ~$h$~, ~$L'(h,k/2+1/4) \neq 0$ and
the real number 
\begin{equation*}
\exp\left(\frac{L^{'}(h,k/2+1/4)}{L(h,{k}/{2}+{1}/{4}}+\Psi(k/2+1/4)\right) 
\end{equation*} 
is transcendental. Moreover, considering the quotient ~$\frac{L'(h,k/2+1/4)}{L(h,{k}/{2}+{1}/{4})}$~ as a function of ~$k$ one deduce 
that the function ~$\frac{L'(h,k/2+1/4)}{L(h,{k}/{2}+{1}/{4})}+\Psi(k/2+1/4)$~ is independent of ~$k$~ and the function
~$\frac{L'(h,k/2+1/4)}{L(h,{k}/{2}+{1}/{4})}\rightarrow -\infty$~ as ~$k\rightarrow \infty$.
\end{thm}
\begin{prop}\label{prop:1}
If ~$\frac{L'(h_0,k_0/2+1/4)}{L(h_0,{k_0}/{2}+{1}/{4})}$~ is algebraic (resp. transcendental) for some
~$h_0\in E(k_0+1/2,4N)$~ then ~$\frac{L'(h,k/2+1/4)}{L(h,{k}/{2}+{1}/{4})}$~ is algebraic (resp. transcendental) for all
~$h\in E(k+1/2,4N)$~ and for all ~$k\in \mathbb{N}$~ with ~$k\equiv k_0 \pmod{2}$.
\end{prop}

\noindent {\bf Note:} 
In \cite[\S 4]{NK}, the above results are proved for the case ~$N=1$ and was remarked that a similar method will lead to the results for $N$ square-free.
In fact following the same arguments one can get the above results for any positive integer ~$N$. 

The above results (\thmref{thm3} and \propref{prop:1}) are about the properties of the functions in the set ~$E(k+1/2, 4N)$~ 
under the assumption that it is a non-empty set. 
It is to be noted that the main results of our present work guarantees the existence of an element in the set $E(k+1/2, 4N)$. 
In fact, by taking ~$r_0=0$~ and ~$\sigma=k/2+1/4$~ in \thmref{thm1} and \thmref{thm2}, we get the following two corollaries. 

\begin{cor}
There exists a constant $C_1$ such that for any $k>C_1$ and any $N$, the set $E(k+1/2,4N)$ is non-empty. 
\end{cor}
\begin{cor}
There exists a constant $C_2>0$ such that for any $N>C_2$ and any  $k\geq 3$, the set $E(k+1/2,4N)$ is non-empty.
\end{cor}

\noindent {\bf Acknowledgement}: The authors thank the referee for making valuable suggestions.

\end{document}